\documentclass[12pt,a4paper,reqno]{amsart}

\usepackage[T1]{fontenc}
\usepackage{lmodern}
\usepackage[utf8]{inputenc}

\usepackage[
    margin=1.3in
]{geometry}

\usepackage{amsmath,amsthm}
\usepackage{amssymb}
\usepackage{mathrsfs}
\usepackage{graphicx}
\graphicspath{{./}{figures/}}

\usepackage{
    color,
    etoolbox,
    pgffor,
    array,
    xcolor,
    microtype
}

\IfFileExists{siunitx.sty}{\usepackage{siunitx}}{}
\IfFileExists{dsfont.sty}{\usepackage{dsfont}}{}
\IfFileExists{bropd.sty}{\usepackage{bropd}}{}

\makeatletter
\newcommand{\needspace}[1]{%
  \ifhmode \ifinner \@LRmoderr \else \par \fi \fi
  \begingroup
    \@tempskipb \lastskip
    \vskip -\@tempskipb
    \setlength\dimen@{#1}%
    \advance\dimen@\@tempskipb
    \vskip\z@\@plus\dimen@
    \penalty -100
    \vskip\z@\@plus -\dimen@
    \vskip\dimen@
    \penalty 9999%
    \vskip -\dimen@
    \vskip\@tempskipb 
  \endgroup
}
\newcommand\Needspace{%
  \ifhmode \ifinner \@LRmoderr \else \par \fi \fi
  \@tempskipb \lastskip
  \vskip -\@tempskipb
  \penalty-100
  \begingroup
  \@ifstar{\@sneedsp@}{\@needsp@}%
}
\newcommand\@sneedsp@[1]{%
    \setlength\dimen@{#1}%
    \advance\dimen@\@tempskipb
    \dimen@ii\pagegoal \advance\dimen@ii-\pagetotal
    \ifdim \dimen@>\dimen@ii
      \break
    \else
      \vskip\@tempskipb  
    \fi
  \endgroup
}
\newcommand\@needsp@[1]{%
    \setlength\dimen@{#1}%
    \advance\dimen@\@tempskipb
    \dimen@ii\pagegoal \advance\dimen@ii-\pagetotal
    \ifdim \dimen@>\dimen@ii
      \ifdim \dimen@ii>\z@
        \vfil
      \fi
      \break
    \else
      \vskip\@tempskipb  
    \fi
  \endgroup
}
\makeatother

\usepackage[shortlabels]{enumitem}

\setlist{
    font=\upshape,
    itemsep=0.25\baselineskip,
    leftmargin=\parindent,
    parsep=0pt,
    topsep=0.25\baselineskip
}

\definecolor{freeviolet}{RGB}{0,0,255}
\definecolor{contactred}{RGB}{0,0,0}
\definecolor{phasegreen}{RGB}{255,255,255}

\newtheoremstyle{plain}
    {6.5pt}
    {6.5pt}
    {\itshape}
    {}
    {\bfseries}
    {.}
    {.5em}
    {}

\newtheoremstyle{definition}
    {6.5pt}
    {6.5pt}
    {\normalfont}
    {}
    {\bfseries}
    {.}
    {.5em}
    {}

\theoremstyle{plain}

\newtheorem{theorem}{Theorem}[section]
\newtheorem{lemma}[theorem]{Lemma}
\newtheorem{corollary}[theorem]{Corollary}
\newtheorem{proposition}[theorem]{Proposition}

\theoremstyle{definition}

\newtheorem{definition}[theorem]{Definition}

\newtheorem{remark}[theorem]{Remark}

\theoremstyle{plain}

\newtheorem{innercustomgeneric}{\customgenericname}
\providecommand{\customgenericname}{}

\newcommand{\newcustomtheorem}[2]{%
    \newenvironment{#1}[1]
    {%
        \renewcommand\customgenericname{#2}%
        \renewcommand\theinnercustomgeneric{##1}%
        \innercustomgeneric
    }
    {\endinnercustomgeneric}
}

\newcustomtheorem{maintheorem}{Theorem}

\numberwithin{equation}{section}

\usepackage[
    pdfa,
    colorlinks,
    allcolors=blue
]{hyperref}

\makeatletter

\let\keywords\relax
\let\msc\relax
\let\correspondingauthor\relax
\let\abstract\relax

\gdef\@keywords{}
\gdef\@msc{}
\gdef\@correspondingauthor{}

\newcommand{\keywords}[1]{%
    \gdef\@keywords{%
        \textit{Keywords:} \ignorespaces#1\unskip.%
    }%
}

\newcommand{\msc}[2][]{%
    \gdef\@msc{%
        \textit{Mathematics Subject Classification 2020:} \ignorespaces#2\unskip%
        \ifx&#1&\else\ (#1)\fi.%
    }%
}

\newcommand{\correspondingauthor}[2]{%
    \gdef\@correspondingauthor{%
        \textit{Corresponding author:} #1,
        \href{mailto:#2}{#2}.%
    }%
}

\newcommand{\ems@titlefootnote}[1]{%
    \insert\footins{%
        \reset@font
        \footnotesize
        \interlinepenalty\interfootnotelinepenalty
        \splittopskip\footnotesep
        \splitmaxdepth\dp\strutbox
        \floatingpenalty\@MM
        \hsize\columnwidth
        \@parboxrestore
        \color@begingroup
            \noindent
            \strut
            #1\par
        \color@endgroup
    }%
}

\newbox\emsabstractbox

\newenvironment{abstract}{%
    \global\setbox\emsabstractbox\vbox\bgroup
        \normalcolor
        \small
        \noindent
        \strut
        \textbf{Abstract.}\enskip
        \ignorespaces
}{%
        \unskip
        \strut
        \par
    \egroup
}

\renewcommand{\maketitle}{%

    \thispagestyle{empty}

    \begingroup

        \vspace*{26pt}

        \centering

        {\Large\bfseries\boldmath
            \@title
            \par
        }

        \vspace{8mm}

        {\large
            \@author
            \par
        }

        \vspace{5mm}

    \endgroup

    \ifx\@msc\@empty
    \else
        \ems@titlefootnote{\@msc}%
    \fi

    \ifx\@keywords\@empty
    \else
        \ems@titlefootnote{\@keywords}%
    \fi

    \ifx\@correspondingauthor\@empty
    \else
        \ems@titlefootnote{\@correspondingauthor}%
    \fi

    \ifvoid\emsabstractbox
    \else
        \noindent
        \unvbox\emsabstractbox
        \par
    \fi

    \vspace{26pt}
}

\newcounter{ems@authorcount}
\gdef\@author{}
\gdef\ems@authorids{}

\newcommand{\Author}[3]{%
    \stepcounter{ems@authorcount}%
    \listgadd{\ems@authorids}{#1}%
    \begingroup
        \def\givenname##1{\csgdef{auth@#1@first}{##1}}%
        \def\surname##1{\csgdef{auth@#1@last}{##1}}%
        \def\mrid##1{%
            \ifblank{##1}{}{\csgdef{auth@#1@mrid}{##1}}%
        }%
        \def\zblid##1{%
            \ifblank{##1}{}{\csgdef{auth@#1@zblid}{##1}}%
        }%
        \def\orcid##1{%
            \ifblank{##1}{}{\csgdef{auth@#1@orcid}{##1}}%
        }%
        #2%
    \endgroup
    \ifnum\value{ems@authorcount}=1
        \gdef\@author{%
            \csuse{auth@#1@first}\space\csuse{auth@#1@last}%
        }%
    \else
        \g@addto@macro\@author{%
            ,\space\csuse{auth@#1@first}\space\csuse{auth@#1@last}%
        }%
    \fi
}

\def\ems@setmeta#1#2#3{%
    \ifblank{#3}{}{%
        \csgdef{affil@#1@#2}{#3}%
        \csgdef{affil@#1@present}{1}%
    }%
}

\def\ems@makeaffilmacro#1{%
    \expandafter\def\csname #1\endcsname##1{%
        \@ifnextchar\bgroup
            {\csuse{ems@#1@two}{##1}}%
            {\csuse{ems@#1@one}{##1}}%
    }%
    \expandafter\def\csname ems@#1@two\endcsname##1##2{%
        \ems@setmeta{\ems@currentaffil@id @##1}{#1}{##2}%
    }%
    \expandafter\def\csname ems@#1@one\endcsname##1{%
        \ems@setmeta{\ems@currentaffil@id @1}{#1}{##1}%
    }%
}

\newcommand{\Affil}[2]{%
    \begingroup
        \def\ems@currentaffil@id{#1}%
        \ems@makeaffilmacro{department}%
        \ems@makeaffilmacro{organisation}%
        \ems@makeaffilmacro{address}%
        \ems@makeaffilmacro{post}%
        \ems@makeaffilmacro{city}%
        \ems@makeaffilmacro{country}%
        \ems@makeaffilmacro{affemail}%
        #2%
    \endgroup
}

\def\ems@emitfield#1{%
    \ifcsdef{affil@\ems@affkey @#1}{%
        \ems@sep
        \csuse{affil@\ems@affkey @#1}%
        \def\ems@sep{,\space}%
    }{}%
}

\def\ems@printid#1#2#3{%
    \ifcsdef{auth@\ems@aid @#1}{%
        \edef\ems@tempid{\csuse{auth@\ems@aid @#1}}%
        \ems@idsep
        #2~\href{#3\ems@tempid}{\ems@tempid}%
        \def\ems@idsep{\quad}%
    }{}%
}
\def\ems@checkid#1{%
    \ifcsdef{auth@\ems@aid @#1}{\def\ems@hasid{1}}{}%
}

\def\ems@printauthor#1{%
    \begingroup
        \def\ems@aid{#1}%
        \Needspace{5\baselineskip}%
        \noindent
        {\bfseries\csuse{auth@#1@first}\space\csuse{auth@#1@last}}%
        \par\nobreak
        \foreach \j in {1,2,3,4}{%
            \edef\ems@affkey{\ems@aid @\j}%
            \ifcsdef{affil@\ems@affkey @present}{%
                \def\ems@sep{}%
                \ems@emitfield{department}%
                \ems@emitfield{organisation}%
                \ifx\ems@sep\@empty\else\par\nobreak\fi
                \def\ems@sep{}%
                \ems@emitfield{address}%
                \ems@emitfield{post}%
                \ems@emitfield{city}%
                \ems@emitfield{country}%
                \ifx\ems@sep\@empty\else\par\fi
                \ifcsdef{affil@\ems@affkey @affemail}{%
                    \edef\ems@tempemail{%
                        \csuse{affil@\ems@affkey @affemail}%
                    }%
                    \href{mailto:\ems@tempemail}{\ems@tempemail}\par
                }{}%
            }{}%
        }%
        \def\ems@hasid{0}%
        \forcsvlist{\ems@checkid}{zblid,mrid,orcid}%
        \ifnum\ems@hasid=1
            \def\ems@idsep{}%
            \noindent Author IDs:\space
            \ems@printid{zblid}{zbMATH}{https://zbmath.org/authors/}%
            \ems@printid{mrid}{MR}{https://mathscinet.ams.org/mathscinet/MRAuthorID/}%
            \ems@printid{orcid}{ORCID}{https://orcid.org/}%
            \par
        \fi
        \addvspace{13pt}%
    \endgroup
}

\AtEndDocument{%
    \ifnum\value{ems@authorcount}>0
        \par\addvspace{19.5pt}%
        \begingroup
            \small\raggedright
            \forlistloop{\ems@printauthor}{\ems@authorids}%
        \endgroup
    \fi
}

\makeatother

\newcommand{\dx}[1]{%
    \mathop{}\!\mathrm{d}#1
}

\newcommand{\dH}[1]{%
    \mathop{}\!\mathrm{d}\mathcal{H}^{#1}
}

\def\XXint#1#2#3{%

    {%
        \setbox0=\hbox{%
            $#1{#2#3}{\int}$%
        }%

        \vcenter{%
            \hbox{$#2#3$}%
        }%

        \kern-.53\wd0
    }%
}

\newcommand{\R}{\mathbb{R}}

\newcommand{\cH}{\ensuremath{\mathcal H}}

\let\ge\geqslant
\let\le\leqslant

\newcommand{\mres}{%
    \mathbin{%
        \vrule
            height 1.6ex
            depth 0pt
            width 0.13ex
        \vrule
            height 0.13ex
            depth 0pt
            width 1.3ex
    }%
}

\DeclareMathOperator{\Div}{div}
\title{%
    Nonexistence of cusp singularities for gravity waves of extreme form
}

\Author{1}{%

    \givenname{Lili}
    \surname{Du}

    \mrid{}
    \zblid{}
    \orcid{}

}{}

\Author{2}{%

    \givenname{Yang}
    \surname{Pu}

    \mrid{}
    \zblid{}
    \orcid{}

}{}

\Author{3}{%

    \givenname{Chunlei}
    \surname{Yang}

    \mrid{}
    \zblid{}
    \orcid{}

}{}

\Affil{1}{%

    \department{Department of Mathematics}

    \organisation{Sichuan University}

    \post{610065}

    \city{Chengdu}

    \country{China}

    \affemail{dulili@scu.edu.cn}
}

\Affil{2}{%

    \department{School of Mathematics}

    \organisation{Southwestern University of Finance and Economics}

    \post{610031}

    \city{Chengdu}

    \country{China}

    \affemail{puyang@swufe.edu.cn}
}

\Affil{3}{%

    \department{School of Mathematical Sciences}

    \organisation{Shenzhen University}

    \post{518060}

    \city{Shenzhen}

    \country{China}

    \affemail{yangchunlei@szu.edu.cn}
}

\msc{%
    35R35, 35B65, 76B15
}


\correspondingauthor
    {Yang Pu}
    {puyang@swufe.edu.cn}

\hypersetup{pdftitle={Nonexistence of cusp singularities for gravity waves of extreme form},pdfauthor={Lili Du, Yang Pu, Chunlei Yang}}
\begin{document}
\enlargethispage{2pt}
\begin{abstract}
In this paper we revisit the geometric approach developed by
V\u{a}rv\u{a}ruc\u{a} and Weiss for the Stokes conjecture for gravity water
waves near stagnation points, in both the irrotational and rotational settings.
They introduced monotonicity formulas and frequency formulas to classify all
possible singularities of solutions at stagnation points. Two questions were
left open:
\begin{itemize}
    \item whether cusp singularities can be excluded in higher dimensions
        \(n\ge 3\) [\emph{Acta Math.} 206 (2011)];
    \item whether cusp singularities for rotational waves can be excluded under
        the Rayleigh-Taylor condition
        [\emph{Ann.\ Inst.\ H.\ Poincar\'e C Anal.\ Non Lin\'eaire} 29 (2012)].
\end{itemize}
We answer both questions in the affirmative. Our main results
show that cusp singularities do not exist in higher dimensions and cannot occur
for rotational waves even without assuming the Rayleigh-Taylor condition.

Our argument combines a delicate analysis of the weighted volume arising from
the corresponding free boundary problem with a weighted isoperimetric inequality.
\end{abstract}

\maketitle

\section{Introduction}\label{sec.introduction}

The investigation of singularities in gravity water waves lies at the
intersection of partial differential equations and fluid mechanics.
Mathematically, the physical problem can be translated into the study of the
one-phase semilinear Bernoulli-type problem
\begin{equation}\label{eq.problem}
\left\{
\begin{aligned}
\Delta u&=-f(u)\quad&&\text{in }\Omega\cap\{u>0\},\\
|\nabla u|^2&=x_n\quad&&\text{on }\Omega\cap\partial\{u>0\}
\end{aligned}
\right.
\end{equation}
where \(\Omega\subset\R^n\) (\(n\ge2\)) is a bounded open set meeting the hyperplane
\(\{x_n=0\}\) and \(f\) is a real-valued continuous function. In two dimensions,
the unknown \(u\) in \eqref{eq.problem} is the stream function and \(f\) the
vorticity-strength function. When \(f=0\), the flow is irrotational; we refer the reader
to the introductions of \cite{VW11,VW12} for a precise derivation of
\eqref{eq.problem} from the Euler equations in two dimensions.

From the perspective of water waves, we are restricted to \(n=2\) and are
particularly interested in solutions of \eqref{eq.problem} with singularities at
the stagnation points, at which the relative fluid velocity \((u_y,-u_x)\)
vanishes. Here subscripts denote partial derivatives. It follows from the second
equation of \eqref{eq.problem} that such points lie on the intersection
\(\partial\{u>0\}\cap\{x_2=0\}\). On the other hand, from the point of view of
free boundary problems, the Bernoulli free boundary condition \(|\nabla u|^2=x_n\)
leads us to consider the natural rescaling
\begin{equation}\label{eq.rescaling}
u_r(x):=\frac{u(x^0+rx)}{r^{3/2}}\quad\text{ for every }x^0\in\partial\{u>0\}\cap\{x_n=0\}.
\end{equation}
This scaling preserves the Bernoulli free boundary condition.

The singular behavior of the free surface near stagnation points is closely
related to the conjecture proposed by G. Stokes \cite{Sto80}, which predicts a
symmetric corner of angle \(2\pi/3\) at the crest of the wave. When \(f=0\),
this conjecture was first proved independently in \cite{AFT82,Plo02}, under
structural assumptions including the isolation of stagnation points and the
symmetry and monotonicity of the free surface. The rotational case was proved in
\cite{Var09} under the same structural assumptions, together with sign
conditions on the vorticity function \(f\).

\begin{figure}[htbp]
\centering
\includegraphics[width=0.6\linewidth]{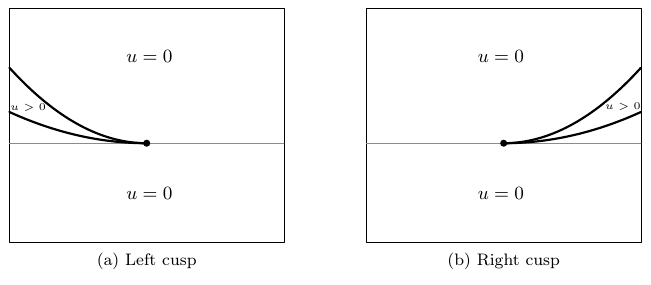}
\caption{The cusp figure.}
\label{cusp_figure}
\end{figure}

In 2011--2012, V\u{a}rv\u{a}ruc\u{a} and Weiss \cite{VW11,VW12} studied the
problem from a geometric point of view, without relying on the structural assumptions that were required in previous works. They established monotonicity formulas
and analyzed the rescaling \eqref{eq.rescaling}. Moreover, under the local growth
assumption
\begin{equation}\label{eq.growth}
|\nabla u|^2\le C_0x_n^+\quad\text{locally in }\Omega\cap\{u>0\},\qquad C_0>0,
\end{equation}
they proved that \(u_r\to u_0\) strongly in \(W_{\rm loc}^{1,2}(\R^n)\) and locally
uniformly in \(\R^n\) along a subsequence as \(r\to0^+\). If
\(u_0\not\equiv 0\), then \(u_0\) is a homogeneous function of degree \(3/2\)
and \(\{u_0>0\}\) is an open sector of angle \(120^\circ\) in two dimensions.
On the other hand, if \(u_0\equiv 0\), they introduced the
\emph{Weiss weighted density}, defined as the limit
\begin{equation}\label{eq.weiss_density}
M_{x^0}(0^+)=\lim_{r\rightarrow0^+}r^{-n-1}\int_{B_r(x^0)}x_n^+\chi_{\{u>0\}}\,\dx{x}.
\end{equation}
Here and below, \(\chi_A\) denotes the characteristic function of a measurable set \(A\).
A stagnation point \(x^0\in\Omega\cap\partial\{u>0\}\cap\{x_n=0\}\) is
\emph{a cusp point} if \(M_{x^0}(0^+)=0\), meaning that the Weiss weighted
density at \(x^0\) vanishes, see Figure \ref{cusp_figure}.

On the other hand, it is called a
\emph{horizontally flat singularity} if
\[
M_{x^0}(0^+)=\int_{B_1}x_n^+\,\dx{x}.
\]
To give a rigorous mathematical proof of the Stokes conjecture, one needs to exclude these two singularities. The nonexistence of horizontally flat singularities has been established for the
irrotational case when the zero set \(\{u=0\}\) has only finitely many
components, and has been proved for the rotational case when the vorticity
satisfies \(|f(z)|\le Cz\) \cite{VW12}; it has also been proved in \cite{DY26a}
under some weaker structural assumptions on \(f\). This topic is beyond the scope of
the present paper, and we focus only on cusp singularities.

The exclusion of cusp singularities for the irrotational case (\(f=0\) in
\eqref{eq.problem}) relies on two key ingredients \cite[Lemma 4.4]{VW11}: the
strong Bernstein estimate \(|\nabla u|^2\le x_2^+\), that is, \(C_0=1\) in
\eqref{eq.growth}, and the local \(C^{2,\alpha}\)
regularity of the topological free boundary
\(\partial\{u>0\}\cap\Omega\cap\{x_2>\tau\}\) outside a set of
\(\mathcal H^1\)-measure zero for every \(\tau>0\)\footnote{See Definition 3.2 in \cite{VW11}
for the definition of a weak solution.}. The key inequality in this method is
\[
 \Delta u\ge\sqrt{x_2}\,\cH^1\mres\partial_{\rm red}\{u>0\}\mbox{\footnotemark},
\]
and the geometric property is that, for cusp points, a non-empty portion of
\(\partial_{\rm red}\{u>0\}\) touches the stagnation point \(x^0\).
\footnotetext{\(\cH^1\mres\partial_{\rm red}\{u>0\}\) denotes the
restriction of the one-dimensional Hausdorff measure to the reduced
boundary of \(\{u>0\}\).}
Then the
authors asked the following question \cite[Page 380]{VW11}:
\begin{quotation}
    \textit{In higher dimensions it is not so clear whether cusps can be excluded. Of course that does not happen for Lipschitz free boundaries.}
\end{quotation}
The nonexistence of cusp singularities for Lipschitz free boundaries was proved
in \cite[Lemma 4.6]{VW11} under the growth condition \eqref{eq.growth}.
Recently, McCurdy \cite[Theorem 1.2]{McC25} excluded cusp singularities in two
dimensions under the weaker growth condition \eqref{eq.growth} by introducing a
height function for the free boundary near the cusp singularity. However,
the exclusion of cusp singularities in higher dimensions remained open prior to the present work.
Let us stress that the primary aim of this paper is to rule out cusp singularities in higher dimensions, thereby answering the open question raised in \cite[Page 380]{VW11}.

The Stokes conjecture for waves with non-zero vorticity was considered by
V\u{a}rv\u{a}ruc\u{a} and Weiss in \cite{VW12}. There they extended their
geometric approach from \cite{VW11} to the semilinear Bernoulli-type problem
\eqref{eq.problem} (\(f\not\equiv 0\) in
\eqref{eq.problem}). Their analysis showed that the limiting problem reduces to
the irrotational case. However, the analysis of the cusp singularities
introduces new difficulties due to the non-zero right-hand side, which causes
additional challenges compared to the irrotational case. Even under the strong
Bernstein estimate \(|\nabla u|^2\le x_2^+\), the method of
\cite[Lemma 4.4]{VW11} breaks down here. Consequently, they conjectured that the
cusp alternative is also not possible under the strong Rayleigh-Taylor condition
\cite[Page 872]{VW12}.
\begin{quotation}
\textit{Unfortunately, we do not have the Rayleigh-Taylor condition}
\begin{equation}\label{eq.strong}
|\nabla u|^2+2F(u)-x_n^+\le 0,\quad\text{ where }\quad F(s)=\int_0^sf(t)\,\dx{t}.
\end{equation}
\textit{in the case with non-zero vorticity, and the method of \cite{VW11}
breaks down here. Still we conjecture that the cusps in case (iii) are not
possible when assuming the Rayleigh-Taylor condition.}
\end{quotation}
Let us remark that there are two sources of difficulty for this rotational
problem. One difficulty is that the combination of vorticity and gravity makes it hard to obtain
the estimate \eqref{eq.strong}. On top of that, we cannot follow the argument of
\cite[Lemma 4.4]{VW11} to exclude the cusp singularity even if we assume
\eqref{eq.strong}. Consequently, the second goal of this paper is to exclude cusp singularities in the two-dimensional case, settling this open question.

In this paper, we prove that cusp singularities cannot occur by showing that the
Weiss weighted density \eqref{eq.weiss_density} has a positive lower bound:
\[
M_{x^0}(0^+)\ge c(n,C_0)>0
\]
for every stagnation point \(x^0\in\partial\{u>0\}\cap\Omega\cap\{x_n=0\}\),
where \(c(n,C_0)\) is a constant depending only on \(n\) and \(C_0\). Under the growth condition \eqref{eq.growth}, the lower bound holds
in every dimension \(n\ge2\) for variational solutions in the sense
of \cite{VW11,VW12}. When the flow is
rotational, we assume only that the vorticity function \(f\) is continuous.
Moreover, our results do not rely on any structural assumptions including the
isolation of stagnation points and the monotonicity and symmetry of the free
boundary. Furthermore, our results show that the Rayleigh-Taylor condition is
in fact not necessary for the nonexistence of cusp singularities.

The paper is organized as follows. In Section~\ref{sec.main}, we state the main
results. In
Section~\ref{sec.preliminaries}, we fix the notation and collect the Weiss-type
monotonicity formula and the weighted isoperimetric inequality. In
Section~\ref{sec.proof}, we give the proof of Theorem \ref{thm.main}.

\section{Main results}\label{sec.main}

Before we state the main results, we recall the notion of variational solutions
of \eqref{eq.problem}. Assume that the nonlinearity \(f\) in \eqref{eq.problem}
is a continuous function with primitive \(F(s)=\int_0^s f(t)\,\dx{t}\).

\begin{definition}[Variational solutions, see Definition 3.1 in \cite{VW12}]\label{def.variational}
We define \(u\in W_{\rm loc}^{1,2}(\Omega)\) to be a \emph{variational solution}
of \eqref{eq.problem} if \(u\in C^0(\Omega)\cap C^2(\Omega\cap\{u>0\})\),
\(u\ge0\) in \(\Omega\) and \(u\equiv 0\) in \(\Omega\cap\{x_n\le0\}\), and the
first variation with respect to domain variations of the functional
\[
J(v):=\int_\Omega\bigl(|\nabla v|^2-2F(v)+x_n\chi_{\{v>0\}}\bigr)\,\dx{x}
\]
vanishes at \(v=u\), i.e.,
\begin{equation}\label{eq.stationarity}
0=\int_\Omega\bigl((|\nabla u|^2-2F(u))\Div\phi
-2\nabla uD\phi\nabla u\bigr)\,\dx{x}+\int_\Omega(x_n\Div\phi+\phi_n)\chi_{\{u>0\}}\,\dx{x}
\end{equation}
for every \(\phi\in C_c^1(\Omega;\R^n)\).
\end{definition}

\begin{remark}
Since our results are completely local, we do not specify any boundary
conditions on \(\partial\Omega\).
\end{remark}

\begin{remark}
If \(f\equiv 0\), then the notion of variational solutions coincides with
the one in \cite[Definition 3.1]{VW11} for the irrotational case.
\end{remark}

We write \(x=(x',x_n)\in\mathbb R^{n-1}\times\mathbb R\), \(x_n^+=\max\{x_n,0\}\) and denote
by \(B_r(x^0)\) the open ball of radius \(r\) centered at \(x^0\).

Throughout the paper, we assume the local growth condition
\eqref{eq.growth}. Fix \(x^0\in\Omega\cap\{x_n=0\}\) and
\(B_R(x^0)\Subset\Omega\). With \(C_0\) denoting the growth constant
on this ball, set
\begin{equation}\label{eq.M0}
M_0=\max\left\{f^+(s):0\le s\le\frac23\sqrt{C_0}R^{3/2}\right\},
\end{equation}
where \(f^+(s)=\max\{f(s),0\}\).
\begin{remark}
Let \(u\) be a variational solution of \eqref{eq.problem}. Since \(u\in C^2(B_R(x^0)\cap\{u>0\})\), the domain-stationarity identity
\eqref{eq.stationarity} holds for test functions supported in the positive
phase, so we can deduce that
\[(|\nabla u|^2-2F(u))I-2\nabla u\otimes\nabla u\]
is divergence free in the positive phase. Here \(\nabla u\otimes\nabla u\) is the symmetric matrix with entries \((\nabla u\otimes\nabla u)_{ij}=(\partial_i u)(\partial_j u)\). The divergence of the displayed tensor equals \(-2(\Delta u+f(u))\nabla u\), hence
\[
(\Delta u+f(u))\nabla u=0\quad\hbox{in }B_R(x^0)\cap\{u>0\}.
\]
If \(\Delta u+f(u)\) were nonzero at some positive point, this identity would
make \(u\) constant on the connected component of
\(B_R(x^0)\cap\{u>0\}\cap\{\Delta u+f(u)\ne0\}\) containing it. Write
\(u=a>0\) on this component. Then \(\Delta u=0\) and
\(\Delta u+f(u)=f(a)\ne0\) there. At every limit point of the component inside
\(B_R(x^0)\), continuity of \(u\) gives \(u=a>0\), and continuity of its second
derivatives in the positive phase gives \(\Delta u=0\). Such a point therefore
still belongs to \(\{u>0\}\cap\{\Delta u+f(u)\ne0\}\) and to the same connected
component. That component would then be both open and closed in \(B_R(x^0)\),
and since \(B_R(x^0)\) is connected it would equal \(B_R(x^0)\), contradicting
\(u=0\) in the lower half-ball. Therefore
\begin{equation}\label{proof.interior.equation}
\Delta u=-f(u)\quad\hbox{in }B_R(x^0)\cap\{u>0\}.
\end{equation}
\end{remark}

We now
define the set of stagnation points and the set of cusp singularities.

\begin{definition}[Stagnation points]\label{def.stagnation_points}
Let \(u\) be a variational solution of \eqref{eq.problem}.
The set of stagnation points is defined by
\[
 \mathcal S^u
 :=\{x\in\Omega\cap\partial\{u>0\}:x_n=0\}.
\]
\end{definition}

\begin{remark}
It should be noted that Definition \ref{def.stagnation_points} introduces a generalized
notion of stagnation points, which corresponds to the physical concept of
stagnation points only when \(n=2\).
\end{remark}

\begin{definition}[Cusp singularities]\label{def.cusp}
For a variational solution \(u\) of \eqref{eq.problem}, the set of
cusp singularities is
\[
 \mathcal C^u=\{x^0\in\mathcal S^u:M_{x^0}(0^+)=0\},
\]
where \(M_{x^0}(0^+)\) is the Weiss weighted density in \eqref{eq.weiss_density}.
\end{definition}

\begin{remark}
The existence and finiteness of the Weiss weighted density are proved in
\cite[Lemma 4.4]{VW12} for the rotational case and in \cite[Lemma 4.2]{VW11} for
the irrotational case, both under the growth condition \eqref{eq.growth}.
Moreover, if \(x^0\in \mathcal C^u\) is a cusp singularity, then
\[
u_r\to u_0\equiv 0\quad\text{ as }r\to 0^+,
\]
strongly in \(W_{\rm loc}^{1,2}(\R^n)\) and locally uniformly in \(\R^n\), where
\(u_r\) is the natural scaling defined in \eqref{eq.rescaling}.
\end{remark}

We now state our main results.

\begin{maintheorem}{I}\label{thm.main}
Let \(n\ge2\), let \(u\) be a variational solution of \eqref{eq.problem} and
\(B_R(x^0)\Subset\Omega\) be a ball with \(x^0\in\mathcal S^u\). Assume that the
growth condition \eqref{eq.growth} holds a.e. in \(B_R(x^0)\). Then there exists a constant \(c(n,C_0)>0\) depending only on  \(n\) and \(C_0\) such that
\begin{equation}\label{eq.quantitative}
r^{-n-1}\int_{B_r(x^0)}x_n^+\chi_{\{u>0\}}\,\dx{x}
\ge c(n,C_0)\left(1-\frac23\sqrt{C_0}M_0\sqrt r\right)^{2n+1}
\end{equation}
for every \(0<r<R\) such that \(\sqrt{C_0}M_0\sqrt r<1\), with \(M_0\) defined in
\eqref{eq.M0}. Consequently,
\[
\mathcal C^u=\varnothing.
\]
\end{maintheorem}

\begin{remark}
It follows from the definition of the Weiss weighted density in \eqref{eq.weiss_density} and the estimate \eqref{eq.quantitative} that
\[
M_{x^0}(0^+)\ge c(n,C_0)>0.
\]
In addition, the asymptotic lower bound is independent of the vorticity. When
\(f\le0\) on the relevant interval, \(M_0=0\) and \eqref{eq.quantitative} holds
for every \(0<r<R\).
\end{remark}

As a direct corollary of Theorem~\ref{thm.main}, we obtain the following
vanishing result for cusp singularities.

\begin{corollary}\label{cor.vanishing}
Let \(n\ge 2\), let \(u\) be a variational solution of \eqref{eq.problem}, and
let \(f\in C([0,\infty))\). Assume that the growth condition \eqref{eq.growth}
holds near \(x^0\in\Omega\cap\{x_n=0\}\). Then \(M_{x^0}(0^+)=0\) implies that \(u\equiv 0\) in
some open \(n\)-dimensional ball containing \(x^0\).
\end{corollary}

\begin{proof}
Since \(x_n^0=0\), we have \(u(x^0)=0\). If \(u\) did not vanish in any
neighborhood of \(x^0\), every ball centered at \(x^0\) would contain a point of
\(\{u>0\}\). Continuity would then give \(x^0\in\partial\{u>0\}\).
Theorem~\ref{thm.main} would imply \(M_{x^0}(0^+)\ge c(n,C_0)>0\), contradicting
the hypothesis. Hence \(u\) vanishes in some ball centered at \(x^0\).
\end{proof}

\begin{remark}
The main novelty of our approach is a new isoperimetric argument that excludes cusp singularities in arbitrary dimensions and for general continuous vorticity, without the Rayleigh–Taylor condition. Indeed, through Theorem \ref{thm.main} and Corollary \ref{cor.vanishing}, we answer two open questions raised in \cite[Page 380]{VW11} and \cite[Remark 4.7]{VW12}, respectively:
\begin{itemize}
    \item The nonexistence of cusp singularities in higher dimensions when \(f\equiv 0\);
    \item The nonexistence of cusp singularities for general continuous \(f\) without assuming the Rayleigh-Taylor condition.
\end{itemize}
\end{remark}

\begin{remark}
The continuity of \(f\) is not necessary; it can be weakened to \(f\in C((0,\infty))\cap L^1((0,\infty))\) with \(f^+\) bounded near \(0\). The bound \eqref{eq.quantitative} still holds and rules out cusps.
\end{remark}

\begin{remark}In the irrotational case the pressure inside the fluid is superharmonic, which yields the Rayleigh–Taylor condition \(|\nabla u|^2\le x_n^+\) on smooth parts of the free boundary. For general vorticity, however, the pressure need not be superharmonic, and our argument uses neither this property nor the Rayleigh–Taylor condition. Our approach combines the domain variation identity \eqref{eq.stationarity} with a weighted isoperimetric inequality to derive directly a positive lower bound for the Weiss weighted density, a quantity that provides an important geometric characterization of cusp singularities.
\end{remark}

\section{Preliminaries}\label{sec.preliminaries}

In the rest of the paper, we write \(B_r\) for \(B_r(0)\). The outward unit normal to
\(\partial B_r\) is \(\nu=x/r\), and we denote the volume of the unit
ball by \(\omega_n=|B_1|\).

For a measurable set \(A\), we denote its characteristic function by
\(\chi_A\) and its distributional derivative by \(\nabla\chi_A\).
On any open set where \(A\) has locally finite perimeter,
\(|\nabla\chi_A|\) denotes the total variation measure of \(\nabla\chi_A\). If \(A\) has
a smooth boundary with outward unit normal \(\nu_A\), then
\[
 \nabla\chi_A=-\nu_A\,\mathcal H^{n-1}\mres\partial A,
\]
where \(\mathcal H^{n-1}\) denotes the \((n-1)\)-dimensional Hausdorff
measure. We use the standard local definition of \(BV\) and follow the
notation of \cite{EG15,AFP00}. Unless otherwise specified, perimeters are taken relative to the
upper half-space \(\{x_n>0\}\).

Once the weighted perimeter is known to be finite, we can define
the weighted volume and the interior weighted perimeter of the positive
phase in \(B_r\) by
\begin{equation}\label{eq.VP}
\begin{aligned}
 V(r)&=\int_{B_r\cap\{x_n>0\}}\sqrt{x_n}\,\chi_{\{u>0\}}\,\dx{x},\\
 P(r)&=\int_{B_r\cap\{x_n>0\}}\sqrt{x_n}\,
       \dx{|\nabla\chi_{\{u>0\}}|}.
\end{aligned}
\end{equation}
As \(u\in W^{1,2}_{\mathrm{loc}}(\Omega)\),
\cite[Theorem~4.4(iv)]{EG15} gives
\begin{equation}\label{proof.zero.gradient}
 \nabla u=0
 \quad\text{almost everywhere on }\{u=0\}\cap B_R(x^0).
\end{equation}
The growth bound \eqref{eq.growth} therefore holds almost everywhere
in \(B_R(x^0)\), and \(u\) is locally Lipschitz there. After translating \(x^0\) to zero, the vertical segment joining
\((x',0)\) to \((x',x_n)\in B_R\cap\{x_n>0\}\) stays in the ball. Thus, for almost every \(x'\), we have
\[
u(x',x_n)=\int_0^{x_n}\partial_nu(x',s)\,\dx{s}
\le\sqrt{C_0}\int_0^{x_n}\sqrt s\,\dx{s}.
\]
Continuity extends the result to every point, so
\begin{equation}\label{eq.vertical.growth}
0\le u(x)\le\frac23\sqrt{C_0}(x_n^+)^{3/2}\qquad x\in B_R.
\end{equation}
It follows from \eqref{eq.M0} that \(f^+(u)\le M_0\) in \(B_R\). The same bound therefore holds in every
concentric ball \(B_r\) with \(0<r<R\).

We now record the monotonicity formula for variational solutions \(u\) of
\eqref{eq.problem}.

\begin{proposition}\label{prop.weiss}
Let \(n\ge 2\), let \(u\) be a variational solution of \eqref{eq.problem}, and let \(f\in C([0,\infty))\). Assume that the growth condition
\eqref{eq.growth} holds in \(B_R(x^0)\Subset\Omega\) with \(x_n^0=0\). Define
\[
\begin{aligned}
M_{x^0}(r)={}&r^{-n-1}\int_{B_r(x^0)}
\bigl(|\nabla u|^2-uf(u)+x_n^+\chi_{\{u>0\}}\bigr)\,\dx{x}\\
&-\frac32r^{-n-2}\int_{\partial B_r(x^0)}u^2\,\dH{n-1},
\qquad 0<r<R.
\end{aligned}
\]
Then \(M_{x^0}\) is locally absolutely continuous on \((0,R)\) and
\begin{equation}\label{eq.weiss.derivative}
M_{x^0}'(r)=2r^{-n-1}\int_{\partial B_r(x^0)}
\left(\nabla u\cdot\nu-\frac{3u}{2r}\right)^2\,\dH{n-1}
+r^{-n-2}K(r)
\end{equation}
for a.e. \(r\in(0,R)\), where
\begin{equation}\label{eq.K}
K(r)=r\int_{\partial B_r(x^0)}(2F(u)-uf(u))\,\dH{n-1}+\int_{B_r(x^0)}\bigl((n-2)uf(u)-2nF(u)\bigr)\,\dx{x}.
\end{equation}
Moreover, we have
\begin{enumerate}[label=(\roman*)]
    \item \(K(r)\) satisfies \(|r^{-n-2}K(r)|\le Cr^{-1/2}\) for small \(r>0\).
    \item The limit \(M_{x^0}(0^+)\) exists, is finite, and satisfies \eqref{eq.weiss_density}.
    \item Every sequence \(r_j\rightarrow0^+\) has a subsequence along which
    \(u_{r_j}\) converges locally uniformly and strongly in
    \(W^{1,2}_{\mathrm{loc}}(\R^n)\) to a \(3/2\)-homogeneous function, harmonic in
    its positive phase and zero on \(\{x_n\le0\}\).
\end{enumerate}
\end{proposition}

\begin{proof}
The formulas \eqref{eq.weiss.derivative} and \eqref{eq.K} were proved in \cite[Theorem~3.4]{VW12}, the error estimate for \(K(r)\) is
\cite[Remark~4.2]{VW12}. The density formula follows from \cite[Lemma~4.4(i),(iv)]{VW12}.
For (iii), the growth estimates give a locally uniformly convergent,
weakly \(W^{1,2}_{\mathrm{loc}}\)-convergent subsequence, and
\cite[Lemma~4.4(ii),(iii)]{VW12} gives strong convergence and
\(3/2\)-homogeneity of the limit \(u_0\).
On compact subsets of \(\{u_0>0\}\), the right-hand side of the
rescaled equation tends uniformly to zero, so \(u_0\) is harmonic there.
\end{proof}

\begin{remark}
For variational solutions, the convergence \(u_r\to u_0\equiv0\) as \(r\to0^+\) does not
determine the limiting behavior of \(\chi_{\{u_r>0\}}\). By Proposition~\ref{prop.weiss}(ii)
and the rescaling \eqref{eq.rescaling},
\[
 M_{x^0}(0^+)
 =\lim_{r\to0^+}\int_{B_1}x_n^+\chi_{\{u_r>0\}}\,\dx{x}.
\]
Thus the Weiss density records the limiting weighted volume of the
rescaled positive phase even in the case of a zero blow-up limit.
\end{remark}

In the following lemma, we will introduce the sharp isoperimetric inequality for sets in the upper half-space with weight \(\sqrt{x_n}\).
This lemma is proved in \cite[Theorem~1.4]{CR13}, while the version for measurable sets can be founded in \cite[Theorem~1.3]{CRS16}.

\begin{lemma}\label{lem.isoperimetric}
Let \(A\subset\{x_n>0\}\) be bounded and measurable with finite weighted
perimeter. Then
\begin{equation}\label{eq.isoperimetric}
\begin{aligned}
 \int_{\{x_n>0\}}\sqrt{x_n}\,\dx{|\nabla\chi_A|}
 \ge{}&\left(n+\frac12\right)
 \left(\int_{B_1\cap\{x_n>0\}}\sqrt{x_n}\,\dx{x}\right)^{\frac{2}{2n+1}}\\
 &\times\left(\int_A\sqrt{x_n}\,\dx{x}\right)^{\frac{2n-1}{2n+1}}.
\end{aligned}
\end{equation}
\end{lemma}

In addition to the weighted isoperimetric inequality, we will use the following comparison between the ordinary volume and the weighted perimeter. It will allow us to control the volume term arising from the vorticity and absorb it at sufficiently small scales.

If \(A\subset B_r\cap\{x_n>0\}\) is a measurable set of finite perimeter and \(\overline A\subset\{x_n>0\}\), then
\begin{equation}\label{eq.volume.comparison}
|A|\le\sqrt r\int_{\{x_n>0\}}\sqrt{x_n}\,\dx{|\nabla\chi_A|}.
\end{equation}

To prove this, choose
\[
\eta\in C_c^\infty(\{x_n>0\}),
\qquad
\eta\equiv1
\quad\text{in a neighborhood of }\overline A,
\]
and set
\[
\phi(x)=\eta(x)x_ne_n,
\]
where \(e_n\) denotes the \(n\)-th coordinate unit vector. Then
\[
\operatorname{div}\phi=1
\quad\text{in a neighborhood of }\overline A.
\]
The Gauss--Green formula therefore gives
\[
\begin{aligned}
|A|
=&-\int_{\{x_n>0\}}\phi\cdot\dx{(\nabla\chi_A)}\\
=&-\int_{\{x_n>0\}}\eta x_n\,\dx{(\partial_n\chi_A)}\\
=&-\int_{\{x_n>0\}} x_n\,\dx{(\partial_n\chi_A)}.
\end{aligned}\]
Here the last equality follows from \(\eta=1\) on \(\operatorname{supp}|\nabla\chi_A|\subset\overline A\).
Since \(A\subset B_r\cap\{x_n>0\}\) has finite perimeter and \(0<x_n\le r\) on this support, we obtain
\[
|A|\le\int_{\{x_n>0\}}x_n\,\dx{|\nabla\chi_A|}\le\sqrt r\int_{\{x_n>0\}}\sqrt{x_n}\,\dx{|\nabla\chi_A|}.
\]

We shall also need some preliminary estimates that recover the boundary measure and the weighted perimeter from the domain variation identity \eqref{eq.stationarity} satisfied by the variational solution \(u\).

\begin{lemma}\label{lem.measure}
Let \(n\ge2\), let \(f\in C([0,\infty))\), and let \(u\) be a variational
solution of \eqref{eq.problem} in \(\Omega\). Assume that
\(B_R\Subset\Omega\) and
\[
|\nabla u|^2\le C_0x_n^+\quad\text{a.e. in }B_R
\]
for some \(C_0>0\).
Then \[\mu=\Delta u+f(u)\chi_{\{u>0\}}\,\dx{x}\] is a
nonnegative Radon measure on \(B_R\) whose support is contained in \(B_R\cap\partial\{u>0\}\). Moreover, for every \(r\) with \(0<r<R\),
\[\mu(B_r)\le C(n)\sqrt{C_0}\,r^{n-1/2}+C(n)M_0r^n.\]
\end{lemma}

\begin{proof}
Since \(f(u)\chi_{\{u>0\}}\) is locally integrable in \(B_R\),
we define the distribution \(\mu\) on \(B_R\) by
\begin{equation}\label{proof.mu}
 \langle\mu,\zeta\rangle
 =-\int_{B_R}\nabla u\cdot\nabla\zeta\,\dx{x}
  +\int_{B_R}f(u)\chi_{\{u>0\}}\zeta\,\dx{x},
 \qquad \zeta\in C_c^\infty(B_R).
\end{equation}
By \eqref{proof.mu} and the interior equation
\eqref{proof.interior.equation}, \(\langle\mu,\zeta\rangle=0\) whenever
\(\zeta\) is supported either in \(B_R\cap\{u>0\}\) or in
\(B_R\cap\operatorname{int}\{u=0\}\). It follows that
\[
 \operatorname{supp}\mu\subset B_R\cap\partial\{u>0\},
\]
where the support is taken relative to \(B_R\).
To prove that \(\mu\) is nonnegative,
choose a smooth nondecreasing function \(\eta:\R\to[0,1]\) such that
\[
\eta(s)=
\begin{cases}
0, & s\le 1,\\
1, & s\ge 2.
\end{cases}
\]
For \(\varepsilon>0\),
\(\eta(u/\varepsilon)\zeta\) is compactly supported in
\(B_R\cap\{u>0\}\), we can use it as a test function in
\eqref{proof.interior.equation}. Integrating by parts gives
\begin{equation}\label{proof.scalar.cutoff}
\begin{aligned}
 \frac1\varepsilon\int_{B_R}
       \zeta\eta'(u/\varepsilon)|\nabla u|^2\,\dx{x}
 ={}&-\int_{B_R}\eta(u/\varepsilon)
                        \nabla u\cdot\nabla\zeta\,\dx{x}\\
 &+\int_{B_R}\zeta\eta(u/\varepsilon)f(u)\,\dx{x}.
\end{aligned}
\end{equation}
Since \(\eta'(u/\varepsilon)\ge0\), we can define a nonnegative Radon measure
\[
 \mu_\varepsilon
 =\frac{1}{\varepsilon}\eta'(u/\varepsilon)|\nabla u|^2\,\dx{x}.
\]
Letting \(\varepsilon\to0^+\), the cutoff \(\eta(u/\varepsilon)\) converges
pointwise to \(\chi_{\{u>0\}}\). By \eqref{proof.zero.gradient}, it is easy to check that
\[\chi_{\{u>0\}}\nabla u=\nabla u \text{ almost everywhere in }B_R.\]
Applying dominated convergence to both terms on the right-hand side of
\eqref{proof.scalar.cutoff} and using \eqref{proof.mu}, we obtain
\[
 \lim_{\varepsilon\to0^+}\int_{B_R}\zeta\,\dx{\mu_\varepsilon}
 =\langle\mu,\zeta\rangle
 \qquad\text{for every }\zeta\in C_c^\infty(B_R).
\]
For \(\zeta\ge0\), the integrals on the left are nonnegative, so
\(\langle\mu,\zeta\rangle\ge0\). Thus \(\mu\) is represented by a
nonnegative Radon measure on \(B_R\).

We now extend the convergence to continuous compactly supported
test functions. For each compact set \(K\Subset B_R\), choose a
nonnegative \(\xi\in C_c^\infty(B_R)\) equal to one on \(K\).
By the convergence just proved for smooth test functions,
\[
\limsup_{\varepsilon\to0^+}\mu_\varepsilon(K)
 \le \lim_{\varepsilon\to0^+}\int_{B_R}\xi\,\dx{\mu_\varepsilon}
 =\int_{B_R}\xi\,\dx{\mu}<\infty.
\]
Hence \(\mu_\varepsilon(K)\) is uniformly bounded for all sufficiently
small \(\varepsilon\). Every function in \(C_c(B_R)\) can be approximated uniformly by
smooth functions whose supports lie in a common compact subset of
\(B_R\). Thus,
\begin{equation}\label{proof.mu.convergence}
 \mu_\varepsilon\rightharpoonup\mu
 \quad\text{locally as Radon measures in }B_R
\end{equation}
Fix \(0<r<R\) and choose \(r<r_1<\min\{2r,R\}\).
Take \(\zeta_r\in C_c^\infty(B_{r_1})\) such that
\[
 \zeta_r=1 \quad\text{on }B_r,\qquad
 0\le\zeta_r\le1,\qquad
|\nabla\zeta_r|\le \frac{C}{r_1-r}.
\]
Since \(r_1<2r\), we have
\[
 |B_{r_1}\setminus B_r|\le C(n)r^{n-1}(r_1-r).
\]
Using the nonnegativity of \(\mu\), the gradient bound, and
\(f^+(u)\le M_0\) in $B_R$, we obtain
\[
\begin{aligned}
 \mu(B_r)
 &\le \int_{B_{r_1}}\zeta_r\,d\mu\\
 &=-\int_{B_{r_1}}\nabla u\cdot\nabla\zeta_r\,\dx{x}
   +\int_{B_{r_1}}f(u)\chi_{\{u>0\}}\zeta_r\,\dx{x}\\
 &\le \frac{C\sqrt{C_0r_1}}{r_1-r}\,|B_{r_1}\setminus B_r|
       +M_0|B_{r_1}|\\
 &\le C(n)\sqrt{C_0}\,r^{n-1/2}+C(n)M_0r^n.
\end{aligned}
\]
\end{proof}

\begin{lemma}\label{lem.perimeter}
Suppose that the hypotheses of Lemma~\ref{lem.measure} hold.
Let \(\mu\) be the nonnegative Radon measure from Lemma~\ref{lem.measure}. Then \[\chi_{\{u>0\}}\in BV_{\mathrm{loc}}(B_R\cap\{x_n>0\}),\]
and
\begin{equation}\label{proof.weighted.perimeter}
 \sqrt{x_n}\,|\nabla\chi_{\{u>0\}}|
\le\sqrt{C_0}\,\mu
      \quad\hbox{in }B_R\cap\{x_n>0\}.
\end{equation}
\end{lemma}

\begin{proof}
Let \(\phi\in C_c^1(B_R\cap\{x_n>0\};\R^n)\) and let \(\eta\) be
the same cutoff function as in the proof of Lemma~\ref{lem.measure}.
It is easy to check that \(\eta(u/\varepsilon)=0\) wherever
\(u\le\varepsilon\) for \(\varepsilon>0\). Moreover, we have
\[
 \operatorname{supp}\bigl(\eta(u/\varepsilon)\phi\bigr)
 \subset \operatorname{supp}\phi\cap\{u\ge\varepsilon\}
 \Subset B_R\cap\{u>0\}.
\]
Multiplying the \(\Delta u+f(u)=0\) by
\(2\eta(u/\varepsilon)\nabla u\cdot\phi\) and integrating over
\(B_R\), we have
\[
0={}2\int_{B_R}
     \eta(u/\varepsilon)\,\Delta u\,(\nabla u\cdot\phi)\,\dx{x}
   +2\int_{B_R}
     \eta(u/\varepsilon)\,f(u)(\nabla u\cdot\phi)\,\dx{x}.
\]
Integrating by parts, the first term becomes
\[
\begin{aligned}
2\int_{B_R}\eta(u/\varepsilon) \Delta u\,(\nabla u\cdot\phi)\,\dx{x}
={}&-\frac{2}{\varepsilon}\int_{B_R}
  \eta'(u/\varepsilon)|\nabla u|^2(\nabla u\cdot\phi)\,\dx{x}\\
&-2\int_{B_R}\eta(u/\varepsilon)\,
  \nabla u\cdot\nabla(\nabla u\cdot\phi)\,\dx{x}.
\end{aligned}
\]
Note that
\[
 2\nabla u\cdot\nabla(\nabla u\cdot\phi)
 =\nabla(|\nabla u|^2)\cdot\phi
  +2\nabla u D\phi\nabla u.
\]
Integrating by parts in the term containing \(\nabla(|\nabla u|^2)\)
gives
\[
\begin{aligned}
 -\int_{B_R}\eta(u/\varepsilon)\nabla(|\nabla u|^2)\cdot\phi\,\dx{x}=&\int_{B_R}\eta(u/\varepsilon)|\nabla u|^2\Div\phi\,\dx{x}\\
 &+\frac{1}{\varepsilon}\int_{B_R}
   \eta'(u/\varepsilon)|\nabla u|^2(\nabla u\cdot\phi)\,\dx{x}.
\end{aligned}
\]
Hence
\[
\begin{aligned}
 2\int_{B_R}\eta(u/\varepsilon)\Delta u\,(\nabla u\cdot\phi)\,\dx{x}
 ={}&\int_{B_R}\eta(u/\varepsilon)
 \bigl(|\nabla u|^2\Div\phi-2\nabla u D\phi\nabla u\bigr)\,\dx{x}\\
 &-\frac{1}{\varepsilon}\int_{B_R}
 \eta'(u/\varepsilon)|\nabla u|^2(\nabla u\cdot\phi)\,\dx{x}.
\end{aligned}
\]
Using \(f(u)\nabla u=\nabla F(u)\), we also obtain
\[
\begin{aligned}
 2\int_{B_R} \eta(u/\varepsilon)f(u)(\nabla u\cdot\phi)\,\dx{x}
 &=-2\int_{B_R} \eta(u/\varepsilon)F(u)\Div\phi\,\dx{x}\\
 &\quad-\frac{2}{\varepsilon}\int_{B_R}
   \eta'(u/\varepsilon)F(u)(\nabla u\cdot\phi)\,\dx{x}.
\end{aligned}
\]
Adding the two integration by parts identities above,
we arrive at
\begin{equation}\label{proof.stress.cutoff}
\begin{aligned}
 &\int_{B_R}\eta(u/\varepsilon)
   \bigl((|\nabla u|^2-2F(u))\Div\phi
               -2\nabla uD\phi\nabla u\bigr)\,\dx{x}\\
 &\qquad=\frac{1}{\varepsilon}\int_{B_R}\eta'(u/\varepsilon)
          \bigl(|\nabla u|^2+2F(u)\bigr)
          \nabla u\cdot\phi\,\dx{x}\\
 &\qquad=\int_{B_R}\nabla u\cdot\phi\,\dx{\mu_\varepsilon}
       +\frac{2}{\varepsilon}\int_{B_R}\eta'(u/\varepsilon)
                    F(u)\nabla u\cdot\phi\,\dx{x}.
\end{aligned}
\end{equation}
The same dominated convergence argument as in the
proof of Lemma~\ref{lem.measure} shows that the left-hand side of
\eqref{proof.stress.cutoff} converges to
\[
 \int_{B_R}
 \bigl((|\nabla u|^2-2F(u))\Div\phi
       -2\nabla u D\phi\nabla u\bigr)\,\dx{x}
\]
as $\varepsilon\rightarrow0^{+}$.
For the first term on the right, the gradient bound \eqref{eq.growth} and
\eqref{proof.mu.convergence} give
\[
\begin{aligned}
 \limsup_{\varepsilon\rightarrow0^{+}}
 \left|\int_{B_R}\nabla u\cdot\phi\,\dx{\mu_\varepsilon}\right|
 &\le \sqrt{C_0}\lim_{\varepsilon\rightarrow0^{+}}
       \int_{B_R}\sqrt{x_n^+}\,|\phi|\,\dx{\mu_\varepsilon}\\
 &=\sqrt{C_0}\int_{B_R}\sqrt{x_n^+}\,|\phi|\,\dx{\mu}.
\end{aligned}
\]
Wherever \(\eta'(u/\varepsilon)\ne0\), we have
\(\varepsilon<u<2\varepsilon\), so \(|F(u)|\le Cu\le 2C\varepsilon\)
for sufficiently small \(\varepsilon>0\). The second term on the right
therefore satisfies
\[
 \left|\frac{2}{\varepsilon}\int_{B_R}\eta'(u/\varepsilon)
       F(u)\nabla u\cdot\phi\,\dx{x}\right|
 \le C\int_{\{\varepsilon<u<2\varepsilon\}\cap\operatorname{supp}\phi}
       |\nabla u|\,|\phi|\,\dx{x}.
\]
The right-hand side tends to zero by dominated convergence as
\(\varepsilon\to0^+\).
Passing to the limit
in that identity and using \eqref{eq.stationarity}, we obtain
\[
\begin{aligned}
 \lim_{\varepsilon\to0^+}
 \int_{B_R}\nabla u\cdot\phi\,\dx{\mu_\varepsilon}
 &=
 \int_{B_R}
 \bigl((|\nabla u|^2-2F(u))\Div\phi
       -2\nabla u D\phi\nabla u\bigr)\,\dx{x}\\
 &=-\int_{B_R}\chi_{\{u>0\}}\Div(x_n\phi)\,\dx{x},
\end{aligned}
\]
where the last equality uses
\(\Div(x_n\phi)=x_n\Div\phi+\phi_n\).
Taking absolute values in this identity and using the preceding estimate, we conclude that
\begin{equation}\label{proof.weighted.duality}
 \left|\int_{B_R}\chi_{\{u>0\}}\Div(x_n\phi)\,\dx{x}\right|
 \le \sqrt{C_0}\int_{B_R}\sqrt{x_n^+}\,|\phi|\,\dx{\mu}.
\end{equation}
Fix a nonempty open set \(U\Subset B_R\cap\{x_n>0\}\). Since \(x_n\) is bounded away
from zero on \(U\), the field \(\phi/x_n\), extended by zero
outside \(U\), is admissible in \eqref{proof.weighted.duality}.
Applying that inequality to \(\phi/x_n\) gives
\begin{equation}\label{proof.BV.duality}
 \left|\int_U \chi_{\{u>0\}}\Div\phi\,\dx{x}\right|
 \le \sqrt{C_0}\int_U \frac{|\phi|}{\sqrt{x_n}}\,\dx{\mu}.
\end{equation}
Let \(\delta=\min_{\overline U}x_n>0\). Taking the supremum over all such \(\phi\) with \(\|\phi\|_\infty\le1\), we obtain
\[
 |\nabla\chi_{\{u>0\}}|(U)
 \le \sqrt{C_0}\int_U\frac{1}{\sqrt{x_n}}\,\dx{\mu}
 \le \sqrt{\frac{C_0}{\delta}}\,\mu(\overline U)<\infty.
\]
Since \(U\) is arbitrary, the characterization of \(BV\) and the
regularity of Radon measures give
\[
 \chi_{\{u>0\}}\in BV_{\mathrm{loc}}(B_R\cap\{x_n>0\})
\]
and
\[
 |\nabla\chi_{\{u>0\}}|
 \le \frac{\sqrt{C_0}}{\sqrt{x_n}}\,\mu
 \quad\text{in }B_R\cap\{x_n>0\}.
\]

Multiplying by \(\sqrt{x_n}\) proves \eqref{proof.weighted.perimeter}.
By the definition of \(P(r)\) and
\eqref{proof.weighted.perimeter}, for \(0<r<R\), we have
\begin{equation}\label{proof.P.mu}
\begin{aligned}
 P(r)
 &=\int_{B_r\cap\{x_n>0\}}\sqrt{x_n}\,
             \dx{|\nabla\chi_{\{u>0\}}|} \\
 &\le \sqrt{C_0}\,\mu\bigl(B_r\cap\{x_n>0\}\bigr)\\
 &\le \sqrt{C_0}\,\mu(B_r)\\
 &<\infty.
\end{aligned}
\end{equation}
The second inequality follows from the monotonicity of \(\mu\), and
\(\mu(B_r)\le\mu(\overline{B_r})<\infty\) because \(\overline{B_r}\)
is a compact subset of \(B_R\).
\end{proof}

\section{Proof of Theorem~\ref{thm.main}}\label{sec.proof}

In this section, we prove Theorem I. We first use the flux across
\(\partial B_r\) to estimate the interior weighted perimeter \(P(r)\)
in terms of \(V'(r)\) and a vorticity error controlled by the ordinary
volume \(|B_r\cap\{u>0\}|\). To absorb this error at small radii, we
truncate the positive phase above the plane \(x_n=\varepsilon\) and
apply the volume comparison \eqref{eq.volume.comparison}. Applying the
 inequality
\eqref{eq.isoperimetric} to the same truncated sets then leads to a
differential inequality for \(V\). Integration gives a lower bound
for \(V(r)\) of order \(r^{n+1/2}\), from which the positive lower
bound for the Weiss weighted density follows by classical Cauchy--Schwarz inequality.
We begin with the following estimate.

\begin{lemma}\label{lem.flux}
Let \(n\ge2\), let \(f\in C([0,\infty))\), and let \(u\) be a variational
solution of \eqref{eq.problem} in \(\Omega\). Assume that
\(B_R\Subset\Omega\) and
\[
|\nabla u|^2\le C_0x_n^+\quad\text{a.e. in }B_R
\]
for some \(C_0>0\). Let \(\mu\) be the measure given as in Lemma~\ref{lem.measure},
and let
\(V\) and \(P\) be as in \eqref{eq.VP}. Then
\[
 \mu(B_r)
 =\int_{\partial B_r}\nabla u\cdot\nu\,\dH{n-1}
  +\int_{B_r}f(u)\chi_{\{u>0\}}\,\dx{x},
\]
and
\[
 P(r)\le C_0V'(r)+\sqrt{C_0}M_0|B_r\cap\{u>0\}|,
\]
hold for almost every \(r\in(0,R)\),
\end{lemma}

\begin{proof}
The weighted volume \(V\) is locally absolutely continuous on \((0,R)\),
and the coarea formula gives
\begin{equation}\label{proof.V.derivative}
 V'(r)=\int_{\partial B_r}\sqrt{x_n^+}\chi_{\{u>0\}}\,\dH{n-1}
                \quad\hbox{for almost every }r\in(0,R).
\end{equation}
We next relate \(\mu(B_r)\) to the outward flux of \(\nabla u\)
across \(\partial B_r\). Fix \(0<r<R\) and \(0<h<R-r\), define the radial cut-off function
\(\eta_h\) by
\[
\eta_h(x)=
\begin{cases}
1, & |x|\le r,\\[2mm]
\dfrac{r+h-|x|}{h}, & r<|x|<r+h,\\[2mm]
0, & |x|\ge r+h.
\end{cases}
\]
By smooth approximation, \(\eta_h\) is admissible in the distributional
identity \eqref{proof.mu} for \(\mu\). Testing with \(\eta_h\) yields
\[
 \int_{B_R}\eta_h\,\dx{\mu}
 ={}\frac1h\int_{B_{r+h}\setminus B_r}
                   \nabla u\cdot\frac{x}{|x|}\,\dx{x}+\int_{B_R}\eta_h f(u)\chi_{\{u>0\}}\,\dx{x}.
\]
We pass to the limit as \(h\rightarrow0^+\), beginning with the left-hand
side. Dominated convergence implies that
\[
\lim_{h\rightarrow0^+}\int_{B_R}\eta_h\,\dx{\mu}
=\mu(\overline{B_r})
=\mu(B_r)
\]
whenever \(\mu(\partial B_r)=0\), which holds for almost every
\(r\in(0,R)\).

For the first term on the right-hand side, the coarea formula gives
\[
 \frac1h\int_{B_{r+h}\setminus B_r}
 \nabla u\cdot\frac{x}{|x|}\,\dx{x}
 =
 \frac1h\int_r^{r+h}
 \int_{\partial B_s}\nabla u\cdot\nu\,\dH{n-1}\,\dx{s}.
\]
The spherical flux is locally integrable in \(s\). The Lebesgue
differentiation theorem therefore shows that this expression converges to
\[
\int_{\partial B_r}\nabla u\cdot\nu\,\dH{n-1}
\quad\text{for almost every }r\in(0,R).
\]

For the remaining term, the boundedness of \(f(u)\) in \(B_R\) and
dominated convergence give
\[
\lim_{h\rightarrow0^+}
\int_{B_R}\eta_h f(u)\chi_{\{u>0\}}\,\dx{x}
=
\int_{B_r}f(u)\chi_{\{u>0\}}\,\dx{x}.
\]
It follows that
\begin{equation}\label{proof.spherical.flux}
 \mu(B_r)=\int_{\partial B_r}\nabla u\cdot\nu\,\dH{n-1}
 +\int_{B_r}f(u)\chi_{\{u>0\}}\,\dx{x}
\end{equation}
for almost every \(r\in(0,R)\).

By Fubini's theorem in polar coordinates and \eqref{proof.zero.gradient},
the growth bound implies that, for almost every \(r\),
\[
\begin{aligned}
 \int_{\partial B_r}\nabla u\cdot\nu\,\dH{n-1}
 &\le\int_{\partial B_r}|\nabla u|\chi_{\{u>0\}}\,\dH{n-1}\\
 &\le\sqrt{C_0}\int_{\partial B_r}
 \sqrt{x_n^+}\chi_{\{u>0\}}\,\dH{n-1}\\
 &=\sqrt{C_0}V'(r).
\end{aligned}
\]
Moreover, by the definition of \(M_0\),
\[
 \int_{B_r}f(u)\chi_{\{u>0\}}\,\dx{x}
 \le M_0|B_r\cap\{u>0\}|.
\]
Inserting these estimates into \eqref{proof.spherical.flux} and using
\eqref{proof.P.mu}, we conclude that
\begin{equation}\label{proof.flux.perimeter}
 P(r)\le C_0V'(r)+\sqrt{C_0}M_0|B_r\cap\{u>0\}|
\end{equation}
for almost every \(r\in(0,R)\).
\end{proof}

\begin{proposition}\label{prop.truncation}
Suppose that the hypotheses of Lemma~\ref{lem.measure} hold.
Let \(V\) and \(P\) be defined by \eqref{eq.VP}.
For almost every \(r\in(0,R)\),
\[ |B_r\cap\{u>0\}|\le\sqrt r\,[P(r)+V'(r)], \]
and
\[ P(r)+V'(r)\ge\left(n+\frac12\right)\left(\int_{B_1\cap\{x_n>0\}}\sqrt{x_n}\,\dx{x}\right)^{\frac{2}{2n+1}}V(r)^{\frac{2n-1}{2n+1}}. \]
\end{proposition}

\begin{proof}
We first describe the idea. To apply the volume comparison \eqref{eq.volume.comparison}, we
truncate \(B_r\cap\{u>0\}\) above the plane
\(\{x_n=\varepsilon\}\) to produce a bounded set whose
closure lies in the open upper half-space. Once the truncated set is seen to have finite perimeter, both the volume comparison
and the weighted isoperimetric inequality \eqref{eq.isoperimetric}
apply. For almost every \(r\in(0,R)\), as \(\varepsilon\to0^+\),
the interior part of the weighted perimeter tends to \(P(r)\) and
the spherical part tends to \(V'(r)\). The horizontal part tends to zero
because the weight equals \(\sqrt{\varepsilon}\) on
\(B_r\cap\{x_n=\varepsilon\}\), whose \((n-1)\)-dimensional measure
is bounded uniformly in \(\varepsilon\). Thus, the weighted perimeter converges to \(P(r)+V'(r)\), and passing
to the limit in the two inequalities gives the desired estimates.
The proof will be divided into three steps.

\medskip
\noindent\emph{Step 1: Selection of radius.}
Our first aim is to choose a radius \(r\in(0,R)\) such that the trace of
\(\chi_{\{u>0\}}\) taken from inside \(B_r\) equals
\(\chi_{\{u>0\}}(x)\) for \(\mathcal H^{n-1}\)-almost every
\(x\in\partial B_r\cap\{x_n>0\}\).

To see that almost every radius has this property, let \(N\) be the set of points \(x\in B_R\cap\{x_n>0\}\) at which
the Lebesgue value of \(\chi_{\{u>0\}}\) either does not exist or
differs from \(\chi_{\{u>0\}}(x)\). The Lebesgue differentiation
theorem and the coarea formula for \(x\mapsto |x|\) give
\[
 \int_0^R \mathcal H^{n-1}(N\cap\partial B_s)\,\dx{s}=|N|=0.
\]
By the preceding identity and \eqref{proof.V.derivative}, both of the
following equalities hold for almost every \(r\in(0,R)\):
\[
 \mathcal H^{n-1}(N\cap\partial B_r)=0,
 \qquad
 V'(r)=\int_{\partial B_r}\sqrt{x_n^+}\,
              \chi_{\{u>0\}}\,\dH{n-1}.
\]
Fix such a radius \(r\) for the remainder of the proof.
For every \(x\in(\partial B_r\cap\{x_n>0\})\setminus N\),
the ball \(B_r\) has density \(1/2\) at \(x\), which gives
\[
 \frac{|B_\rho(x)|}{|B_\rho(x)\cap B_r|}
 \longrightarrow 2
 \qquad(\rho\rightarrow0^{+}).\]
It follows that
\[
\begin{aligned}
 &\quad\frac{1}{|B_\rho(x)\cap B_r|}
 \int_{B_\rho(x)\cap B_r}
 \bigl|\chi_{\{u>0\}}(y)-\chi_{\{u>0\}}(x)\bigr|\,\dx{y}\\
 &\le
 \frac{|B_\rho(x)|}{|B_\rho(x)\cap B_r|}
 \frac{1}{|B_\rho(x)|}\int_{B_\rho(x)}
 \bigl|\chi_{\{u>0\}}(y)-\chi_{\{u>0\}}(x)\bigr|\,\dx{y}\\
 &\longrightarrow 0
\end{aligned}
\]
as \(\rho\rightarrow0^{+}\). Since \(\mathcal H^{n-1}(N\cap\partial B_r)=0\), this proves the
asserted trace identity for the chosen radius \(r\).

\medskip
\noindent\emph{Step 2. Truncation and decomposition of the weighted perimeter.}
For \(0<\varepsilon<r\), set
\[
 A_\varepsilon=B_r\cap\{u>0\}\cap\{x_n>\varepsilon\}.
\]
Since \(\sqrt{x_n}\ge\sqrt{\varepsilon}\) in
\(B_r\cap\{x_n>\varepsilon\}\), we have
\[
\begin{aligned}
 |\nabla\chi_{\{u>0\}}|\bigl(B_r\cap\{x_n>\varepsilon\}\bigr)
 &\le \frac{1}{\sqrt{\varepsilon}}
       \int_{B_r\cap\{x_n>\varepsilon\}}
       \sqrt{x_n}\,\dx{|\nabla\chi_{\{u>0\}}|}\\
 &\le \frac{P(r)}{\sqrt{\varepsilon}}<\infty.
\end{aligned}
\]
The cap part \(B_r\cap\{x_n>\varepsilon\}\) is a bounded Lipschitz domain. Extending \(\chi_{\{u>0\}}\) from the cap by zero gives
\(\chi_{A_\varepsilon}\). The \(BV\) zero-extension formula yields
\[
|\nabla\chi_{A_\varepsilon}|(\mathbb R^n)
=
|\nabla\chi_{\{u>0\}}|(B_r\cap\{x_n>\varepsilon\})
+
\int_{\partial(B_r\cap\{x_n>\varepsilon\})}
\bigl|\operatorname{Tr}\chi_{\{u>0\}}\bigr|
\,d\mathcal H^{n-1},
\]
where \(\operatorname{Tr}\) denotes the interior trace from \(B_r\cap\{x_n>\varepsilon\}\). Notice that the interior trace of
\(\chi_{\{u>0\}}\) on its boundary takes values in \([0,1]\) almost
everywhere, then
\[
 |\nabla\chi_{A_\varepsilon}|(\mathbb R^n)
 \le \frac{P(r)}{\sqrt{\varepsilon}}
   +\mathcal H^{n-1}\!\left(\partial\bigl(B_r\cap\{x_n>\varepsilon\}\bigr)\right)
 <\infty.
\]
Thus, for every \(0<\varepsilon<r\),
\[
 \chi_{A_\varepsilon}\in BV(\mathbb R^n),
 \qquad
 \overline{A_\varepsilon}\subset\{x_n\ge\varepsilon\}.
\]
Both the weighted isoperimetric inequality \eqref{eq.isoperimetric}
and the volume comparison \eqref{eq.volume.comparison} therefore
apply to \(A_\varepsilon\).
We now employ the \(BV\) Gauss--Green formula to \(\chi_{\{u>0\}}\) in the cap
\(B_r\cap\{x_n>\varepsilon\}\): for every \(\phi\in C_c^1(\R^n;\R^n)\),
\[
\begin{aligned}
\big\langle\nabla\chi_{A_\varepsilon},\phi\big\rangle
&=-\int_{A_\varepsilon}\operatorname{div}\phi\,dx\\
&=\int_{B_r\cap\{x_n>\varepsilon\}}
       \phi\cdot d(\nabla\chi_{\{u>0\}})\\
&\quad-\int_{\partial B_r\cap\{x_n>\varepsilon\}}
       \chi_{\{u>0\}}\,\phi\cdot\nu\,d\mathcal H^{n-1}\\
&\quad+\int_{B_r\cap\{x_n=\varepsilon\}}
       \operatorname{Tr}^{+}\chi_{\{u>0\}}\,\phi_n\,d\mathcal H^{n-1}.
\end{aligned}
\]
Here \(\operatorname{Tr}^{+}\) denotes the trace from above on the
horizontal part. The three measures on the right are supported on the cap
\(B_r\cap\{x_n>\varepsilon\}\), the spherical part
\(\partial B_r\cap\{x_n>\varepsilon\}\), and the horizontal part
\(B_r\cap\{x_n=\varepsilon\}\), respectively. Since these sets are
pairwise disjoint, there total variations add. Integrating this identity against \(\sqrt{x_n}\)
gives
\[
\begin{aligned}
 \int_{\{x_n>0\}}\sqrt{x_n}\,\dx{|\nabla\chi_{A_\varepsilon}|}
 ={}&\int_{B_r\cap\{x_n>\varepsilon\}}
       \sqrt{x_n}\,\dx{|\nabla\chi_{\{u>0\}}|}\\
 &+\int_{\partial B_r\cap\{x_n>\varepsilon\}}
       \sqrt{x_n}\,\chi_{\{u>0\}}\,\dH{n-1}\\
 &+\sqrt{\varepsilon}\int_{B_r\cap\{x_n=\varepsilon\}}
       \operatorname{Tr}^{+}\chi_{\{u>0\}}\,\dH{n-1}.
\end{aligned}
\]
Figure~\ref{fig:trunc} illustrates
the truncation in the planar case.
\begin{figure}[htbp]
\centering
\includegraphics{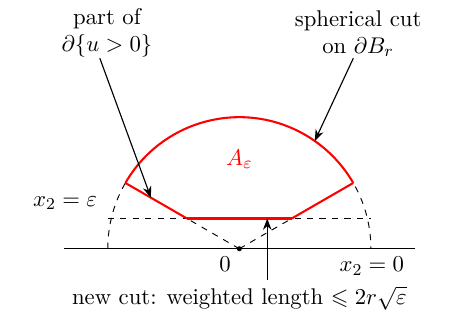}
\caption{The truncated phase \(A_\varepsilon=B_r\cap\{u>0\}\cap\{x_2>\varepsilon\}\).}
\label{fig:trunc}
\end{figure}

\medskip
\noindent\emph{Step 3. Application of the
geometric inequalities.}
As \(\varepsilon\to0^+\), monotone convergence and the identities
\eqref{eq.VP} and \eqref{proof.V.derivative} give
\[
\begin{aligned}
 \int_{B_r\cap\{x_n>\varepsilon\}}
       \sqrt{x_n}\,\dx{|\nabla\chi_{\{u>0\}}|}
 &\rightarrow P(r),\\
 \int_{\partial B_r\cap\{x_n>\varepsilon\}}
       \sqrt{x_n}\,\chi_{\{u>0\}}\,\dH{n-1}
 &\rightarrow
 \int_{\partial B_r}\sqrt{x_n^+}\,\chi_{\{u>0\}}\,\dH{n-1}
 =V'(r).
\end{aligned}
\]
For \(0<\varepsilon<r\), the set \(B_r\cap\{x_n=\varepsilon\}\)
is a disk of radius \(\sqrt{r^2-\varepsilon^2}\). The trace bound
\(0\le\operatorname{Tr}^{+}\chi_{\{u>0\}}\le1\) gives
\[
\begin{aligned}
 0
 &\le \sqrt{\varepsilon}\int_{B_r\cap\{x_n=\varepsilon\}}
       \operatorname{Tr}^{+}\chi_{\{u>0\}}\,\dH{n-1}\\
 &\le \sqrt{\varepsilon}\,
       \mathcal H^{n-1}\bigl(B_r\cap\{x_n=\varepsilon\}\bigr)\\
 &= \sqrt{\varepsilon}\,\omega_{n-1}
       (r^2-\varepsilon^2)^{(n-1)/2}.
\end{aligned}
\]
The horizontal part of the weighted perimeter therefore tends to zero
as \(\varepsilon\to0^+\).
Adding these three limits, we obtain
\[
 \lim_{\varepsilon\rightarrow0^{+}}
 \int_{\{x_n>0\}}\sqrt{x_n}\,\dx{|\nabla\chi_{A_\varepsilon}|}
 =P(r)+V'(r).
\]
Recall that \(u=0\) on \(\{x_n\le0\}\), the sets \(A_\varepsilon\) increase to
\(B_r\cap\{u>0\}\) as \(\varepsilon\to0^+\).
Hence, by monotone convergence,
\[
 |A_\varepsilon|\rightarrow |B_r\cap\{u>0\}|,
 \qquad
 \int_{A_\varepsilon}\sqrt{x_n}\,\dx{x}
 \rightarrow V(r).
\]
By Step~2, the geometric inequalities \eqref{eq.isoperimetric} and
\eqref{eq.volume.comparison} apply to \(A_\varepsilon\) for every
\(0<\varepsilon<r\). Letting \(\varepsilon\rightarrow0^{+}\) in
\eqref{eq.isoperimetric} gives
\begin{equation}\label{proof.radial.isoperimetric}
 P(r)+V'(r)
 \ge
 \left(n+\frac12\right)
 \left(\int_{B_1\cap\{x_n>0\}}\sqrt{x_n}\,\dx{x}\right)^{\frac{2}{2n+1}}
 V(r)^{\frac{2n-1}{2n+1}}.
\end{equation}
Taking the limit in \eqref{eq.volume.comparison} gives
\begin{equation}\label{proof.radial.volume}
 |B_r\cap\{u>0\}|
 \le \sqrt r\,\bigl(P(r)+V'(r)\bigr).
\end{equation}
By Step~1, these inequalities hold for almost every \(r\in(0,R)\).
\end{proof}

\begin{lemma}\label{lem.absorption}
Under the hypotheses of Lemma~\ref{lem.measure}, the estimate
\[
 \bigl(1-\sqrt{C_0}M_0\sqrt r\bigr)\bigl(P(r)+V'(r)\bigr)
 \le (1+C_0)V'(r)
\]
holds for almost every \(r\in(0,R)\).

If, in addition, \(0\in\mathcal S^u\), then for every \(0<r<R\)
such that \(\sqrt{C_0}M_0\sqrt r<1\),
\[
 V(r)\ge
 \left[
   \frac{r}{1+C_0}
   \left(1-\frac23\sqrt{C_0}M_0\sqrt r\right)
 \right]^{n+1/2}
 \int_{B_1\cap\{x_n>0\}}\sqrt{x_n}\,\dx{x}.
\]
In particular, \[V(r)\ge c_1(n,C_0)\,r^{n+1/2},\]
where \(c_1(n,C_0)>0\).
\end{lemma}

\begin{proof}
Insert \eqref{proof.radial.volume} into the flux estimate
\eqref{proof.flux.perimeter}. For almost every radius \(r\in(0,R)\),
\[
 P(r)\le C_0V'(r)
        +\sqrt{C_0}M_0\sqrt r\,[P(r)+V'(r)].
\]
Adding \(V'(r)\) and moving the last term to the left, we obtain
\begin{equation}\label{proof.absorption}
 (1-\sqrt{C_0}M_0\sqrt r)[P(r)+V'(r)]
          \le(1+C_0)V'(r).
\end{equation}
If \(\sqrt{C_0}M_0\sqrt r<1\), the coefficient on the left is positive. We may
therefore substitute the lower bound \eqref{proof.radial.isoperimetric} into
\eqref{proof.absorption} and deduce
\begin{equation}\label{proof.V.inequality}
\begin{aligned}
 V'(r)\ge{}&\frac{n+1/2}{1+C_0}
       (1-\sqrt{C_0}M_0\sqrt r)\\
 &\times\left(\int_{B_1\cap\{x_n>0\}}\sqrt{x_n}\,\dx{x}\right)^{\frac{2}{2n+1}}
       V(r)^{\frac{2n-1}{2n+1}}.
\end{aligned}
\end{equation}
Since \(0\in\partial\{u>0\}\), the set \(B_r\cap\{u>0\}\) is nonempty
for every \(0<r<R\). The weight \(\sqrt{x_n}\) is positive, so
\(V(r)>0\).
Fix \(0<\rho<r<R\) with \(\sqrt{C_0}M_0\sqrt r<1\). On \([\rho,r]\),
\(V\) is absolutely continuous, and its monotonicity gives
\[V(s)\ge V(\rho)>0.\]
The chain rule therefore gives
\[
\begin{aligned}
 \frac{\mathrm d}{\mathrm ds}V(s)^{\frac{2}{2n+1}}
 &=\frac{2}{2n+1}V(s)^{-\frac{2n-1}{2n+1}}V'(s)\\
 &\ge\frac{1-\sqrt{C_0}M_0\sqrt s}{1+C_0}
       \left(\int_{B_1\cap\{x_n>0\}}\sqrt{x_n}\,\dx{x}\right)^{\frac{2}{2n+1}}
\end{aligned}
\]
for almost every \(s\in(\rho,r)\). In the last step, the factor \(n+1/2\) in
\eqref{proof.V.inequality} cancels with \(2/(2n+1)\).

Integrating from \(\rho\) to \(r\) yields
\begin{equation}\label{proof.integrated.positive}
\begin{aligned}
 V(r)^{\frac{2}{2n+1}}-V(\rho)^{\frac{2}{2n+1}}
 \ge{}&\frac1{1+C_0}
 \left(\int_{B_1\cap\{x_n>0\}}\sqrt{x_n}\,\dx{x}\right)^{\frac{2}{2n+1}}\\
 &\times\left(r-\rho
       -\frac23\sqrt{C_0}M_0(r^{3/2}-\rho^{3/2})\right).
\end{aligned}
\end{equation}
Letting \(\rho\to0^+\) in \eqref{proof.integrated.positive} and using the bound
\[
 0<V(\rho)\le\sqrt{\rho}\,|B_\rho\cap\{x_n>0\}|
 =\frac{\omega_n}{2}\rho^{n+1/2},
\]
we obtain
\[
\begin{aligned}
 V(r)^{\frac{2}{2n+1}}
 \ge{}&\frac{r}{1+C_0}
 \left(\int_{B_1\cap\{x_n>0\}}\sqrt{x_n}\,\dx{x}\right)^{\frac{2}{2n+1}}\\
 &\times\left(1-\frac23\sqrt{C_0}M_0\sqrt r\right).
\end{aligned}
\]
The last factor is positive under the radius restriction. Raising both sides to
the power \((2n+1)/2\) gives
\begin{equation}\label{proof.V.lower}
\begin{aligned}
 V(r)\ge{}&\frac{r^{n+1/2}}{(1+C_0)^{n+1/2}}
      \int_{B_1\cap\{x_n>0\}}\sqrt{x_n}\,\dx{x}\\
 &\times\left(1-\frac23\sqrt{C_0}M_0\sqrt r\right)^{n+1/2}.
\end{aligned}
\end{equation}
Since \(\sqrt{C_0}M_0\sqrt r<1\), the last factor in parentheses is greater
than \(1/3\). Thus the asserted lower bound holds with
\[
c_1(n,C_0)=\frac{1}{[3(1+C_0)]^{n+1/2}}
\int_{B_1\cap\{x_n>0\}}\sqrt{x_n}\,\dx{x}>0.
\]
The differential inequality holds for almost every radius, and \(V\) is
absolutely continuous, so its integrated form is valid for every \(0<r<R\)
satisfying \(\sqrt{C_0}M_0\sqrt r<1\).

When \(M_0=0\), the factor \(1-\sqrt{C_0}M_0\sqrt r\) in
\eqref{proof.absorption} equals one. The lower bound for \(V(r)\) then holds for every \(0<r<R\),
and this includes the case \(f\le0\).
When \(M_0>0\), the same factor tends to one as \(r\to0^+\)
and is therefore positive for all sufficiently small \(r\).
\end{proof}
\begin{proof}[Proof of Theorem~\ref{thm.main}]
Fix \(x^0\in\mathcal S^u\) and let \(0<r<R\) satisfy
\(\sqrt{C_0}M_0\sqrt r<1\). After a horizontal translation, we may
assume that \(x^0=0\). To pass from the weight \(\sqrt{x_n}\)
in Lemma~\ref{lem.absorption} to the weight \(x_n\) in
\eqref{eq.quantitative}, we apply the Cauchy--Schwarz inequality on
\(B_r\cap\{u>0\}\), which gives
\[
\begin{aligned}
 V(r)^2
 &=\left(\int_{B_r\cap\{u>0\}}\sqrt{x_n}\,\dx{x}\right)^2\\
 &\le |B_r\cap\{u>0\}|
       \int_{B_r}x_n^+\chi_{\{u>0\}}\,\dx{x}\\
 &\le \frac{\omega_n}{2}r^n
       \int_{B_r}x_n^+\chi_{\{u>0\}}\,\dx{x}.
\end{aligned}
\]
The last inequality follows from \(\{u>0\}\subset\{x_n>0\}\),
since \(B_r\cap\{x_n>0\}\) has volume \(\omega_n r^n/2\).
Combining this inequality with the lower bound \eqref{proof.V.lower}
for \(V(r)\), we obtain
\[
\begin{aligned}
 r^{-n-1}\int_{B_r}x_n^+\chi_{\{u>0\}}\,\dx{x}
 &\ge\frac{2}{\omega_n}\frac{V(r)^2}{r^{2n+1}}\\
 &\ge\frac{2}{\omega_n(1+C_0)^{2n+1}}
       \left(\int_{B_1\cap\{x_n>0\}}\sqrt{x_n}\,\dx{x}\right)^2\\
 &\qquad\times
       \left(1-\frac23\sqrt{C_0}M_0\sqrt r\right)^{2n+1}.
\end{aligned}
\]
This proves \eqref{eq.quantitative} with
\[
 c(n,C_0)=
 \frac{2}{\omega_n(1+C_0)^{2n+1}}
 \left(\int_{B_1\cap\{x_n>0\}}\sqrt{x_n}\,\dx{x}\right)^2>0.
\]
Letting \(r\rightarrow0^+\) yields the positive lower bound
\[
 \liminf_{r\rightarrow0^+}r^{-n-1}
      \int_{B_r}x_n^+\chi_{\{u>0\}}\,\dx{x}
       \ge c(n,C_0)>0.
\]
For continuous vorticity, Proposition~\ref{prop.weiss} gives the existence of the
limiting Weiss weighted density and the identity \eqref{eq.weiss_density}. After translating the center to the origin, we suppress the subscript and write \(M(r)\) for the Weiss functional. Consequently,
\[
 M(0^+)=\lim_{r\rightarrow0^+}r^{-n-1}
      \int_{B_r}x_n^+\chi_{\{u>0\}}\,\dx{x}
       \ge c(n,C_0)>0.
\]
The origin is therefore not a cusp point. Undoing the horizontal translation
proves the conclusion at \(x^0\) and completes the proof of
Theorem~\ref{thm.main}.
\end{proof}

\subsubsection*{Acknowledgement}
\hyphenpenalty=10
\sloppy
L.~Du is supported by National Natural Science Foundation of China under Grants 12125102, 12526202, 12671245, 12631009. Y.~Pu is supported by National Natural Science Foundation of China under Grant 12571226. C.~Yang is supported by National Natural Science Foundation of China under Grant 12601435.

\subsubsection*{AI declarations}
OpenAI’s GPT was used only for grammar, spelling and citation checks, with all edits reviewed to safeguard mathematical reasoning.

\subsubsection*{Conflict of interest}
The authors declare that they have no conflict of interest.

\subsubsection*{Data availability} No data was used in this research.

\end{document}